\documentclass[12 pt]{article}
\usepackage{graphicx}
\usepackage{subfigure}
\usepackage{amsthm}
\newcommand{\hs}{\hskip 0.2cm}
\newcommand{\var}{\varepsilon}

\newcommand{\ha}{\theta}
\newcommand{\va}{\Psi}

\title{\large\bf Invariant Tori of Impulsive Duffing-Type Equations via KAM Technique}
 \author{Lu Chen, Jianhua Shen$^{\thanks{The corresponding author.
{\it E-mail addresses:} chenlu@zju.edu.cn (L. Chen), jianhuashen2013@163.com (J. Shen).}}$\\
\small {\it Department of Mathematics, Hangzhou Normal University}\\
\small {\it Hangzhou,  Zhejiang 310036, PR China}}
\date{}
\begin{document}
 \maketitle

\noindent {\small\bf Abstract} \vskip 0.2cm

{\small  A method via the KAM technique is introduced to study the existence of invariant tori and quasiperiodic solutions for impulsive Duffing-type equations with time period 1. Basing on several planar symplectic homeomorphisms and some estimates of impulsive perturbations under each symplectic homeomorphisms, we prove via the Moser's twist theorem the boundedness (Lagrange stability) and the existence of an invariant circle for the equation with area-preserving impulsive terms. And this invariant circle having any rotation number $\omega>\omega_0$ with some $\omega_0>0$, so we obtain also that the solutions starting from the circle are quasiperiodic with frequencies $\omega$ and 1.

{\small twist theorem}}
 \vskip 0.4cm
 \noindent{\bf 1. Introduction}\vspace{0.2cm}

Consider the second order Duffing-type equation with polynomial potentials
$$\ddot{x}+x^{2n+1}+\sum_{j=0}^{2n} x^{j}p_{j}(t)=0,\hs n\geq1,\eqno (1.1)$$
where $p_j(t)(j=0,1,\cdots,2n)$ are 1-periodic functions. The problem of establishing the boundedness results (Lagrange stability) of solutions for (1.1) was the subject of many important investigations. These studies were driven by the very intricate phenomenon about the longtime behaviour of the nonlinear differential equation $\ddot{x}+f(t,x)=0$ with $f$ being periodic in $t$. For example, there are equations having unbounded solutions but with infinitely many zeroes and with nearby unbounded solutions having randomly prescribed numbers of zeroes and also periodic solutions, see Alekseev [2], Sitnikov [22] and Moser [14,15]. These studies were also driven by the Moser's problem. In fact, Moser [14,15] pointed out that the boundedness of solutions for $\ddot{x}+\alpha x+\beta x^3=p(t)$ could be shown by his twist theorem where $p(t)\in C(S^1),\hs p(t+1)=p(t),\hs S^1=\mathbf{R}/\mathbf{N},\hs \alpha\geq 0,\hs \beta>0$. The first boundedness result is due to Morris [13] who showed that all solutions of $\ddot{x}+2x^3=p(t)$ are bounded, with $p(t)\in C(S^1),\hs p(t+1)=p(t)$. Subsequently, the result was extended to a wider class of systems by Dieckerhoff and Zehnder [4], who showed the boundedness of all solutions and the existence of quasi-periodic solution for (1.1) via the Moser's twist theorem. More results on the studies for the boundedness of (1.1) and some other form via Moser's twist theorem, one may see [7,9,11,12,16,19,24-27].

The underlying idea in [4,24-26] is as follows. By means of transformation theory the equation is outside of a large disc $D$ in the $(x,\dot{x})$-plane transformed into a Hamiltonian equation having the following property. Following the solutions from the section $t=0$ to the section $t=1$ defines a map, the time 1 map $\phi$ of the flow, which is close to a so called twist map in $\mathbf{R^2}\verb|\|D$. By means of the twist-theorem one finds large invariant curves diffeomorphic to circles and surrounding the origin in the $(x,\dot{x})$-plane. Every such curve is the base of a time periodic and under the flow invariant cylinder in the phase space $(x,\dot{x},t)\in \mathbf{R^2}\times \mathbf{R}$, which confines the solutions in its interior and which therefore leads to a bound of these solutions. It turns out that all solutions are Lagrange stable and all the solutions starting at $t=0$ on the invariant curves are quasiperiodic.

It is well known that the invariant curves guaranteed by the twist-theorem lead to an abundance of periodic solutions. For example, by applying the Poincar$\acute{e}$-Birkhoff fixed point theorem to the annuli bounded by two suitable invariant curves one finds fixed points and periodic points of the time 1 map of the flow. They give rise to forced oscillations and subharmonic solutions of the equation. One knows, in addition, that these invariant curves are in the closure of the set of periodic points. However, in [4], by applying the Moser' twist theorem, since the set of invariant curves has infinite Lebesgue measure in $\mathbf{R^2}\setminus D$, the conclusion that the closure of the set of subharmonic solutions of (1.1) is of infinite Lebesgue measure in the phase space $\mathbf{R^2}\times S^1$ could be obtained.

In a KAM technique we need search a symplectic diffeomorphism so that the equation in the $(x,y)$-plane could be transformed into a new Hamiltonian system which is usually described by action and angle-variable form. One of the difficulties we encountered is the determination of an appropriate symplectic transform. Another difficult is the confirm of so called ``small property condition" of Moser's twist theorem about the disturbed part of the new Hamilyonian system. In the processing of the actual equation, the ``intersection property" of the time 1 map is usually guaranteed by the fact that the time 1 map of the Hamiltonian system is area-preserving.

We are concerned in this paper with the existence of invariant tori and quasiperiodic solutions to (1.1) under certain impulsive effect, i.e.,\\
  $$\left\{\begin{array}{ll}
  \ddot{x}+x^{2n+1}+\sum_{i=0}^{2n} x^{i}p_{i}(t)=0, \hs t\neq t_j,\hs n\geq1,\\
  \bigtriangleup x(t_{j})=I_j(x(t_j^-),\dot{x}(t_j^-)),\\
  \bigtriangleup \dot{x}(t_j)=J_j(x(t_j^-),\dot{x}(t_j^-)),\hs j=\pm1,\pm2,\pm3,\cdots,
  \end{array}\right.  \eqno (1.2)$$\\
  where $p_i(t+1)=p_i(t),$  $p_i(t)\in C^{\infty},$  $0\leq t_1<t_2<\cdots<t_{k}<1,$  $\bigtriangleup x(t_{j})=x(t_{j})-x(t_{j}^-),$  $\bigtriangleup \dot{x}(t_j)=\dot{x}(t_j)-\dot{x}(t_j^-)$, and $I_j,J_j:\mathbf{R}\times \mathbf{R}\mapsto \mathbf{R}$ are continuous maps for $j=\pm1,\pm2,\pm3,\cdots$. In addition, assume that the impulsive time is 1-periodic, that is, there exists positive integer $k$ such that $t_{j+k}=t_j+1$ and $I_{j+k}=I_j, J_{j+k}=J_j$ for $j=\pm1,\pm2,\pm3,\cdots$.

  An impulsive differential equation corresponding to a smooth evolution of a dynamics that at certain times changes instantaneously, corresponding to impulses in the smooth system. There are many applications of these equations to mechanical and natural phenomena involving abrupt changes, including in physics, chemistry, biology, control theory, and robotics. We refer the reader to classical monographs [1,8] for the general aspects of impulsive differential equations, [6,17,20] for the existence of periodic solutions via fixed point theory, [5,23] via topological degree theory, and [18] via variational method.

  However, different from the study for Duffing type equations without impulsive terms, there are no results on the invariant tori and quasiperiodic solutions for impulsive Duffing type equations via KAM technique. It is obvious, in a mechanical system, that the periodic impulsive perturbation and periodic external force may both destroy the existence of invariant torus, therefore, studies on the existence of invariant tori for impulsive Duffing type equations are of important significance.

  The main purpose of this paper is to develop a KAM technique to study the existence of invariant tori and quasiperiodic solutions for a Duffing type equation with polynomial potential and nonlinear impulsive terms of the form (1.2) via Moser's twist theorem.

  We note that the method in [4,24-26] now can not be used directly to (1.2), in fact, in the estimation of ``small property condition" of Moser's twist theorem, the calculation of time 1 maps $(\lambda_1,\theta_1)$ of (1.2) in action and angle-variable form and $(\mu_1,\phi_1)$ in a series of symplectic diffeomorphism form will be more complex and difficult because of the complex and difficult calculation of $(\lambda_{t_j^+}, \theta_{t_j^+})$  and $(\mu_{t_j^+}, \phi_{t_j^+})$. Moreover, the determination of regional for the time 1 maps is also depended on those difficult calculation.

   In section 4, we will prove the ``small property condition" by proving that $\frac{\partial^{r+s} \Delta\lambda}{{\partial \lambda^{r}}{\partial \ha^{s}}}=O(\lambda^{-\varepsilon-r}),$
   $\frac{\partial^{r+s} \Delta\ha}{{\partial \lambda^{r}}{\partial \ha^{s}}}=O(\lambda^{-1-\varepsilon-r}),\hs \lambda\rightarrow +\infty,$  if the assumption (i) of Theorem 1.1 holds (see Lemma 4.2), which, in turn, derives $\frac{\partial^{r+s} \Delta\mu}{{\partial \mu^{r}}{\partial \phi^{s}}}=O(\mu^{-\varepsilon-r}),$
   $\frac{\partial^{r+s} \Delta\phi}{{\partial \mu^{r}}{\partial \phi^{s}}}=O(\mu^{-1-\varepsilon-r}),\hs \mu\rightarrow +\infty,$ under a series of symplectic diffeomorphism (see Lemma 4.3).

    Another difficult is the certainty of the ``intersection property" of Poincar$\acute{e}$ time 1 map with impulse effects in case the condition that the jump maps
  $$  \Phi_{j}:(x,y)\mapsto (x,y)+(I_j(x,y),J_j(x,y)),\hs j=1,2,\cdots,k,$$
  is area-preserving homeomorphisms is lacking, because the calculation of Poincar$\acute{e}$ time 1 map (for the intersection-preserving) is very difficult.

  In section 5, we will prove that if the condition (ii) of Theorem 1.1 holds, then
  under a series of symplectic diffeomorphism $\psi_{m}\circ\psi_{m-1}\circ\cdot\cdot\cdot\circ\psi_{1}\circ\va_0$, the time 1 map $P$ of (1.2) is area-preserving, and so the map $P$ has the ``intersection property" (Lemma 5.2).

  The main results in this paper are the following theorems.\\

  \noindent{\bf Theorem 1.1.} {\it Suppose that\\
  {\rm (i)} \hs for any $\var>0$, $j=1,2,\cdots,k$ and any non-negative integers $m_1$ and $m_2$($m_1+m_2\leq5$),
  $$ \lim_{x^{2}+y^2\rightarrow +\infty}\left| \frac{\partial^{m_1+m_2} I_j(x,y)}{{\partial x^{m_1}}{\partial y^{m_2}}}\cdot \left(h_0(x,y)\right)^{\frac{1}{2}+\frac{n+2}{2n+2}\varepsilon+\frac{m_1}{2n+2}+\frac{m_2}{2}}\right|<+\infty,$$
   $$\lim_{x^{2}+y^2\rightarrow +\infty}| \frac{\partial^{m_1+m_2} J_j(x,y)}{{\partial x^{m_1}}{\partial y^{m_2}}}\cdot \left(h_0(x,y)\right)^{\frac{1}{2n+2}+\frac{n+2}{2n+2}\varepsilon+\frac{m_1}{2n+2}+\frac{m_2}{2}}|<+\infty,$$
   where $h_0(x,y)=\frac{1}{2(n+1)}x^{2n+2}+\frac 12 y^2$

  {\rm (ii)}\hs for $j=1,2,\cdots,k$,
  $$\frac{\partial I_j}{\partial x}+\frac{\partial J_j}{\partial y}+\frac{\partial I_j}{\partial x}\cdot \frac{\partial J_j}{\partial y}-\frac{\partial I_j}{\partial y}\cdot \frac{\partial J_j}{\partial x}=0.$$
   Then, all solutions of equation (1.2) are bounded, i.e., for every solution $x(t)$ of equation (1.2),
   $$\sup_{\mathbf{R}}(|x(t)|+|\dot{x}(t)|)<+\infty.$$}

  \noindent {\bf Theorem 1.2.} {\it If all conditions of Theorem 1.1 are satisfied, then there are many irrational numbers $\omega$($\omega>\omega_0$, $\omega_0$ is large), for every $\omega$, there is a quasiperiodic solution of (1.2) having frequencies $(\omega,1)$.}\\

   \noindent {\bf Remark 1.1.} If $I_j, J_j\equiv 0, j=\pm 1, \pm 2, \cdots$, then (1.2) is the second order Duffing-type equation (1.1). Therefore, Theorem 1.1 and 1.2 generalize the corresponding results for ODEs case by Dieckerhoff and Zehnder in [4] to the impulsive setting analogous ones.\\

     \noindent {\bf Remark 1.2.}  It is easy to show that the system of partial differential equations (i.e. condition (ii) of Theorem 1.1 and 1.2) implies that the jump maps
  $$\Phi_{j}:(x,y)\mapsto (x,y)+(I_j(x,y),J_j(x,y)),\hs j=1,2,\cdots  $$
  are area-preserving homeomorphisms. \\

 \noindent {\bf Remark 1.3.} It should be mentioned that the conditions (i) and (ii) in Theorem 1.1 and 1.2 are not easy to verify in $(x,y)$-plane. In section 6, we will give several interesting cases (examples) for the construction of the jump maps $I_j(x,y), J_j(x,y)$ in $(x,y)$-plane by using jump maps in symplectic coordinate form under the symplectic diffeomorphism $\va_0$.\\

\noindent {\bf Remark 1.4.} The equation (1.1) is a Hamiltonian system with the Hamiltonian function
\[H=h_0(x,y)+R(x,y,t)\]
where
\[R(x,y,t)=\sum_{j=0}^{2n}\frac{1}{j+1}p_j(t)x^{j+1}.\]
We point out that $H$ can be generalized to more general function of $(x,y,t)$, not necessarily to be a polynomial of $(x,y)$. Here it is essential to regard $R$ as a relatively small perturbation with respect to $h_0$:
\[\left|\partial_x^{m_1}\partial_{y}^{m_2}\, \left(h_0(x,y)^{-1}\, R(x,y,t)\right) \right|<\infty,\; x^2+y^2\to \infty.\]
See [9] for the detail. Without this relatively small condition imposed on the perturbation $R$, the stability might have been broken. In the present paper, the conditions (i) implies that the impulse is also relatively small with respect to $h_0$. A relatively large impulse might violate the stability, too. So the condition (i) should be reasonable.\\

 \noindent {\bf Remark 1.5.} In the present paper, our interest is in the effect of impulse on the stability. Thus we assume that the coefficients $p_j$'s are $C^\infty$ to simplify the proof. Very recently, Yuan [28] proved that (1.1) is Lagrange stability provided that $p_j(t)$'s are $C^\gamma$-H\"older continuity with $\gamma>1-\frac1n$. By advantage of the idea in Yuan[28], we can prove that the Theorem 1.1 still holds true under condition that all $p_j(t)$'s are $C^\gamma$-H\"older continuity with $\gamma>1-\frac1n$.

 \ \

 The rest of the paper is organized as follows. In Section 2, a standard reduction to the action-angle variables form (see Dieckerhoff and Zehnder [4]) is presented, by which (1.2) is transformed into an impulsive Hamiltonian system.  In Section 3, we will introduce some canonical transformations as in [4]. In section 4, we will prove the ``small property condition"  under a series of symplectic diffeomorphism. In Sections 5, we will prove that the time 1 map $P$ of (1.2) is area-preserving, and so it has the intersection property. Finally, in section 6, we will give some further discussions on the jump maps $I_j$ and $J_j$ for applications.
   \begin{flushleft}
\baselineskip 0.2in \vskip 0.2in{\bf 2. Duffing equation with impulsive terms and its action and angle-variables form
}\setcounter{section}{2} \setcounter{equation}{0}\vskip 0.1in
\end{flushleft}

 In this section, we carry out the standard reduction so that (1.2) will be transformed to its action-angle variables form. Firstly, we recall some basic properties of impulsive differential equations. Consider the initial value problem
 $$
 \left\{
 \begin{array}{ll}
 \dot{u}=F(t,u), \hs t\neq t_j,\\
 \bigtriangleup u_j=L_j(u(t_j^-)), \hs j=\pm 1, \pm 2,\cdots,
  \end{array}
 \right.
 \eqno (2.1)
 $$
 $$  u(t_0^+)=u_0,  \eqno (2.2) $$
where $\bigtriangleup u_j=u(t_j^+)-u(t_j^-), j=\pm 1, \pm 2,\cdots$. Assume that\\

\noindent (i) $F: \mathbf{R}\times \mathbf{R^n}\mapsto \mathbf{R^n}$ is continuous in $(t_j,t_{j+1}]\times \mathbf{R^n}$, locally Lipschitz in the second variable and the limits $\lim_{t\to t_j^+, v\to u}F(t,v),  j=\pm 1, \pm 2,\cdots$ exist;

\noindent (ii) $L_j: \mathbf{R^n}\to \mathbf{R^n},  j=\pm 1, \pm 2,\cdots$, are continuous;

\noindent (iii) $F$ is 1-periodic in the first variable, $0\leq t_1<\cdots<t_k<1, t_{j+k}=t_j+1$ and $L_{j+k}=L_j$ for $ j=\pm 1, \pm 2,\cdots$. \\

We have the following basic lemmas.\\

\noindent {\bf Lemma 2.1.} (See [1,17]) {\it Assume {\rm (i), (ii)} and {\rm (iii)}. Then for any $u_0\in \mathbf{R^n}$, there is a unique solution $u=u(t;t_0,u_0)$ of (2.1) satisfying the initial value condition (2.2). Moreover, $P_t: u_0\mapsto u(t;t_0,u_0)$ is continuous in $u_0$ for $t\neq t_j,  j=\pm 1, \pm 2,\cdots$.}\\

 The second lemma is easily proved.\\

 \noindent {\bf Lemma 2.2.}{\it Assume that the conditions of Lemma 2.1 hold. Then all the solutions of (2.1) exist for $t\in \mathbf{R}$ provided that all the solutions of $\dot{u}=F(t,u)$ exist for $t\in \mathbf{R}$. Moreover, if $\Phi_j: u\mapsto u+L_j(u), j=1,\cdots, k$ are global homeomorphisms of $\mathbf{R^n}$, then the map $P_t$ is a homeomorphism for $t\neq t_j, j=\pm 1, \pm 2,\cdots$. Furthermore, all the solutions of (2.1) have elastic property, that is, for any $b_0>0$, there is $r_{b_0}>0$ such that the inequality $|u_0|\geq r_{b_0}$ implies $|u(t;t_0,u_0)|\geq b_0$, for $t\in (t_0,t_0+1]$.}\\

 Let $t_0=0$ and denote $P_0: u_0\mapsto u(t_1^-;0,u_0), P_j: u_j\mapsto u(t_{j+1}^-;t_j,u_j), j=1,\cdots,k-1, P_k:u_k\mapsto u(1;t_k,u_k)$. Then the Poincar$\acute{e}$ map $P: u_0\mapsto u(1;0,u_0)$ can be expressed by\\
   $$ P=P_k\circ \Phi_k\circ\cdots\circ P_1\circ\Phi_1\circ P_0.   $$

   If $\dot{u}=F(t,u)$ is conservative and $\Phi_j, j=1,\cdots,k$ are global area-preserving homeomorphism, then $P$ is an area-preserving homeomorphism.

 Next, in order to carry out the standard reduction so that (1.2) will be transformed to its action-angle variables form, we write (1.2) as an equivalent system of the form\\
 $$
 \left\{
 \begin{array}{lll}
 \dot{x}=y,\hs t\neq t_j,\\
 \dot{y}=-x^{2n+1}-\sum_{i=0}^{2n}p_i(t)x^i,\hs t\neq t_j;\\
 \bigtriangleup x(t_j)=I_j(x(t_j^-), y(t_j^-)),\\
 \bigtriangleup y(t_j)=J_j(x(t_j^-),y(t_j^-)),\hs j=\pm 1, \pm 2,\cdots.
 \end{array}
 \right.
 \eqno (2.3)
 $$\\
  Denote by $u(t;u_0)=(x(t;u_0),y(t;u_0))$ the solution of (2.3) satisfying the initial condition $u(0^+;u_0)=u_0$, where $u_0=(x_0,y_0)$. Let $P: u_0\mapsto u(1;u_0)$ be the Poincar$\acute{e}$ map of (2.3). Then\\
 $$   P=P_k\circ \Phi_k\circ\cdots\circ P_1\circ\Phi_1\circ P_0,   $$
 where \\
 $$  \Phi_j: u=(x,y)\mapsto (x+I_j(x,y), y+J_j(x,y)), \hs j=1,\cdots, k;   $$
 $$  P_j: u_j\mapsto u_{j+1}^-=u(t_{j+1}^-; u_0),\hs j=0,1,\cdots,k,\hs t_{k+1}=1;  $$
 $$  u_j=\Phi_j(u_j^-),\hs j=1,\cdots,k.    $$\\

 Consider the equation\\
  $$\ddot{x}+x^{2n+1}=0.\eqno (2.4)$$\\
It is clear that we have the field $X_{h}$ of (2.4):\\
$$\left\{\begin{array}{ll}
  \dot{x}=\frac{\partial h(x,y)}{\partial y},\\
  \dot{y}=-\frac{\partial h(x,y)}{\partial x},\\
  h(x,y)=\frac{1}{2}y^2+\frac{1}{2n+2}x^{2n+2}.
  \end{array}\right.  \eqno (2.5)$$\\
It is well known that all nonzero solutions of (2.5) are periodic.

   Suppose that $(x_{0}(t),y_{0}(t))$ is the solution of (2.5) satisfying the initial conditions $(x_{0}(0),y_{0}(0))=(1,0)$. And let $T_{0}$ be its minimal positive period. It follows from (2.5) that $x_{0}(t)$ and $y_{0}(t)$ posses the following qualities:  \\

\noindent (i) $\dot{x}_0(t)=y_0(t)$,\hs $\dot{y}_0(t)=-x_0^{2n+1}(t)$;\\
\noindent (ii) $x_0(t+T_0)=x_0(t)$,\hs $y_0(t+T_0)=y_0(t)$;\\
\noindent (iii) $(n+1)y_0^{2}(t)+x_0^{2n+2}(t)=1$;\\
\noindent (iv) $x_0(-t)=x_0(t)$,\hs $y_0(-t)=-y_0(t)$.\\

  The action-angle variables are now defined by the mapping $\va_0$:$\mathbf{R^+}\times S^1 \longmapsto \mathbf{R^2}\verb|\| \left\{0 \right\}$ , where $(x,y)=\va_0(\lambda,\ha)$ with $\lambda >0$ and $\ha$(mod 1) is given by the formulae:
  $$
  \va_0:
  \left\{\begin{array}{ll}
  x=c^{\alpha}\lambda^{\alpha}x_0(\ha T_0),\\
  y=c^{\beta}\lambda^{\beta}y_0(\ha T_0),
  \end{array}\right.  \eqno (2.6)$$
  where $\alpha=\frac{1}{n+2},\hs \beta=1-\alpha,\hs c=\frac{1}{\alpha T_0}$.
We claim that $\va_0$ is a symplectic diffeomorphism from $\mathbf{R^+}\times S^1$ onto $\mathbf{R^2}\verb|\| \left\{0 \right\}$ (see [4] for the proof). Equation (1.2) has its Hamiltonian function in the $(x,y)$-plane:
  $$h(x,y,t)=\frac{1}{2}y^2+\frac{1}{2n+2}x^{2n+2}+\sum_{j=0}^{2n} \frac{1}{j+1}p_{j}(t)x^{j+1}.\eqno (2.7)$$\\
Under the symplectic transformation $\va_0$, it is therefore transformed into
$$H(\lambda,\ha,t)=\frac{1}{2(n+1)}c^{2\beta}\cdot\lambda^{2\beta}+G(\lambda,\ha,t),\eqno (2.8)$$
where $G$ is of the form:
$$G(\lambda,\ha,t)=\sum_{j=0}^{2n}G_j(\ha,t)\cdot\lambda^{\alpha(j+1)}=\sum_{j=0}^{2n} \frac{p_{j}(t)}{j+1}\cdot c^{\alpha(j+1)}\cdot x_0^{\alpha(j+1)}(\ha T_0)\cdot\lambda^{\alpha(j+1)},$$
with $G_j\in C^{\infty}(T^2)$ and $T^2=R^2/Z^2$.

  By (2.6), we have\\
  $$c\lambda=[x^{2n+2}+(n+1)y^2]^{\frac{n+2}{2n+2}}.\eqno (2.9)$$\\
By (2.3), we have\\
$$\left\{\begin{array}{ll}
  x(t_j)=x(t_j^-)+I_j,\\
  y(t_j)=y(t_j^-)+J_j.
  \end{array}\right.  \eqno (2.10)$$\\
By (2.9) and (2.10), we have\\
$$\lambda(t_j)=\frac{1}{c}[(x(t_j^-)+I_j)^{2n+2}+(n+1)(y(t_j^-)+J_j)^2]^{\frac{n+2}{2n+2}}.\eqno (2.11)$$
Hence, we have
  $$\Delta\lambda(t_j)=\frac{1}{c}[(x(t_j^-)+I_j)^{2n+2}+(n+1)(y(t_j^-)+J_j)^2]^{\frac{n+2}{2n+2}}-\lambda(t_j^-).\eqno (2.12)$$

    Under the symplectic diffeomorphism $\va_0$, equations (1.2) can be transformed into:
  $$\left\{\begin{array}{ll}
  \dot{\ha}=\frac{\partial H}{\partial \lambda},\hs t\neq t_j,\\
  \dot{\lambda}=-\frac{\partial H}{\partial \ha},\hs t\neq t_j,\\
  H(\lambda,\ha,t)=\frac{1}{2(n+1)}c^{2\beta}\cdot\lambda^{2\beta}+G(\lambda,\ha,t),\\
  \bigtriangleup \ha(t_j)=I^{*}_{j}(\lambda(t_j^-),\ha(t_j^-)),\\
  \bigtriangleup \lambda(t_j)=J^{*}_{j}(\lambda(t_j^-),\ha(t_j^-)),\hs j=\pm1,\pm2,\pm3,\cdots,
\end{array}\right.  \eqno (2.13)$$
where $\Delta\ha(t_j)=\ha(t_j)-\ha(t_j^-)$ and $\Delta\lambda(t_j)=\frac{1}{c}[(x(t_j^-)+I_j)^{2n+2}+(n+1)(y(t_j^-)+J_j)^2]^{\frac{n+2}{2n+2}}-\lambda(t_j^-)$.

In the following arguments, we will omit the constant $\frac{1}{2n+2}c^{2\beta}$ without the loss of generality.\\

\noindent {\bf 3. More symplectic transformations} \vskip 0.3cm
We introduce a space of functions $F(r)$ which behave for large $\lambda>0$ as done by Dieckerhoff and Zehnder [4]. Given $r\in \mathbf{R}$ we denote by $F(r)$ the set of $C^{\infty}$ functions in $(\lambda,\ha,t) \in \mathbf{R^+}\times T^2$ which are defined in $\lambda\geq\lambda_{0}$ for some $\lambda_{0}>0$ and for which there is a sequence $\lambda_{jlk}>0$ such that
$$\sup_{\stackrel{\lambda>\lambda_{jlk}}{(\ha,t)\in T^{2}}}(\lambda^{j-r}|(D_{\lambda})^{j}(D_{\ha})^{l}(D_{t})^{k}f(\lambda,\ha,t)|)<\infty.$$
We summarize some properties readily verified from the definition.\\

\noindent{\bf Lemma 3.1.} \\{\it
\noindent {\rm (i)} if  $r_{1}<r_{2}$  then  $F(r_{1})\subset F(r_{2})$ ;\\
\noindent {\rm (ii)} if  $f\in F(r)$  then  $(D_{\lambda})^{j}f\in F(r-j)$ ;\\
\noindent {\rm (iii)} if $f_{1}\in F(r_{1})$  and  $f_{2}\in F(r_{2})$  then  $f_{1}\cdot f_{2}\in F(r_{1}+r_{2})$;\\
\noindent {\rm (iv)} if  $f\in F(r)$  satisfies  $|f(\lambda,\bullet)| \geq c\lambda^{r}$  for  $\lambda>\lambda_{0}$  then  $\frac{1}{f}\in F(-r)$.}\\

  For $f\in F(r)$ we denote the meanvalue over the $\ha$ - variable by $[f]$:
  $$[f](\lambda,t):=\int_{0}^{1}f(\lambda,\ha,t)d\ha.$$

  If $\lambda_0 >0$, then $A_{\lambda_0}\subset R^+\times T^2$ denotes the annulus
  $$A_{\lambda_0}:=\left\{(\lambda,\ha,t)|\lambda\geq\lambda_0,     (\ha,t)\in T^2  \right\} .$$

\noindent{\bf Lemma 3.2} [4]. {\it Let\\
$$ H(\lambda,\theta,t)=\lambda^a+h_1(\lambda,t)+h_2(\lambda,\ha,t),$$\\
with $h_1\in F(c)$ and $h_2\in F(b)$. Assume $a>1$, $b<a$ and $c<a$. Then there is a symplectic diffeomorphism $\psi$ depending periodically on $t$ of the form\\
$$\left\{\begin{array}{ll}
  \lambda=\tilde{\mu}+u(\tilde{\mu},\tilde{\phi},t),\\
  \theta=\tilde{\phi}+v(\tilde{\mu},\tilde{\phi},t),
\end{array}\right.  $$\\
with $u\in F(1-(a-b))$ and $v\in F(-(a-b))$ such that $A_{\tilde{\mu}_+}\subset \psi(A_{\tilde{\mu}_0})\subset A_{\tilde{\mu}_-}$ for some large $\tilde{\mu}_-<\tilde{\mu}_0<\tilde{\mu}_+$. Moreover the transformed Hamiltonian vectorfield $\psi^*(X_{H})=X_{\hat{H}}$ is of the formㄩ
$$\hat{H}(\tilde{\mu},\tilde{\phi},t)=\tilde{\mu}^a+\hat{h}_1(\tilde{\mu},t)+\hat{h}_2(\tilde{\mu},\tilde{\phi},t),$$
where $\hat{h}_1\in F(c_1)$, with $c_1=\max(c,b)$, and where $\hat{h}_2\in F(b_1)$. The constant $b_1$ is smaller than $b$ and is given by
$$
b_1=
\left\{
\begin{array}{ll}
  b-(a-b),\hs  b\geq1,\\
  b-(a-1),\hs  b\leq1.
\end{array}\right.  $$}\\

   We know that system (2.13) has its Hamiltonian function:\\
   $$H(\lambda,\ha,t)=\lambda^a+h_2(\lambda,\ha,t),\eqno (3.1) $$\\
with $h_2\in F(b)$,  $a=\alpha(2n+2)$  and  $b=\alpha(2n+1)$,  where $\alpha=\frac{1}{n+2}$.  Therefore  $1<a<2$  and  $b<a$,  so that the assumptions of Lemma 3.2 are met.\\

\noindent{\bf Lemma 3.3} [4]. {\it Let $H$ be as in (3.1), there is a symplectic diffeomorphism $\psi$, periodic in time t: $\psi=\psi_{m}\circ\psi_{m-1}\circ\cdot\cdot\cdot\circ\psi_{1}$ with $A_{\mu_+}\subset \psi(A_{\mu_0})\subset A_{\mu_-}$ for some large $\mu_-<\mu_0<\mu_+$(m is an integer, $\psi_{i}$ is the symplectic diffeomorphism of Lemma 3.2, $i=1,2,\cdots,m$). Which transforms the Hamiltonian system (3.1) into $\psi^*(X_{H})=X_{\hat{H}}$, where
$$\hat{H}(\mu,\phi,t)=\mu^a+\hat{h}_1(\mu,t)+\hat{h}_2(\mu,\phi,t),\eqno (3.2)$$
with $\hat{h}_1\in F(b)$ and $\hat{h}_2\in F(-\varepsilon_1)$ for $\varepsilon_1>2-a$.}\\

  Under the symplectic diffeomorphism $\psi$ of Lemma 3.3, system (2.13) can be transformed into:
  $$\left\{\begin{array}{ll}
  \dot{\phi}=\frac{\partial \hat{H}}{\partial \mu},\hs t\neq t_j,\\
  \dot{\mu}=-\frac{\partial \hat{H}}{\partial \phi},\hs t\neq t_j,\\
  \hat{H}=\mu^a+\hat{h}_1(\mu,t)+\hat{h}_2(\mu,\phi,t),\\
  \bigtriangleup \phi(t_j)=I^{**}_{j}(\mu(t_j^-),\phi(t_j^-)),\\
  \bigtriangleup \mu(t_j)=J^{**}_{j}(\mu(t_j^-),\phi(t_j^-)),\hs j=\pm1,\pm2,\pm3\cdots\\
\end{array}\right. \eqno (3.3)$$
with $\hat{h}_1\in F(b)$ and $\hat{h}_2\in F(-\varepsilon_1)$ for $\varepsilon_1>2-a$.\\

\noindent {\bf 4. Some estimates} \vskip 0.3cm

\noindent{\bf Lemma 4.1.} {\it  Set $$ \hat{I}_j(\lambda,\ha)=I_j(c^{\alpha}\lambda^{\alpha}x_0(\ha T_0),c^{\beta}\lambda^{\beta}y_0(\ha T_0)),$$
   $$\hat{J}_j(\lambda,\ha)=J_j(c^{\alpha}\lambda^{\alpha}x_0(\ha T_0),c^{\beta}\lambda^{\beta}y_0(\ha T_0)),$$
    then for $j=1,2,\cdots,k$ and any non-negative integers $r$ and $s$($r+s\leq5$),
    $$\frac{\partial^{r+s} \hat{I}_j(\lambda,\ha)}{{\partial \lambda^{r}}{\partial \ha^{s}}}=O(\lambda^{-\frac{n+1}{n+2}-\varepsilon-r}),\hs \lambda\rightarrow +\infty,$$
   $$\frac{\partial^{r+s} \hat{J}_j(\lambda,\ha)}{{\partial \lambda^{r}}{\partial \ha^{s}}}=O(\lambda^{-\frac{1}{n+2}-\varepsilon-r}),\hs \lambda\rightarrow +\infty,$$
if the condition (i) of Theorem 1.1 holds.}\\

\noindent{\bf Proof.} From (2.9) and the condition (i) of Theorem 1.1, we have
$$ \lim_{\lambda\rightarrow +\infty} \frac{\partial^{m_1+m_2} I_j(x,y)}{{\partial x^{m_1}}{\partial y^{m_2}}}=O(\lambda^{-\frac{n+1}{n+2}-\varepsilon-\frac{m_1}{n+2}-\frac{(n+1)m_2}{n+2}}),\eqno (4.1)$$
   $$\lim_{\lambda\rightarrow +\infty} \frac{\partial^{m_1+m_2} J_j(x,y)}{{\partial x^{m_1}}{\partial y^{m_2}}}=O(\lambda^{-\frac{1}{n+2}-\varepsilon-\frac{m_1}{n+2}-\frac{(n+1)m_2}{n+2}}).\eqno (4.2)$$

     {\rm (i)} \hs When $r=0,s=0$, from (4.1) and (4.2), we have
    $$ \hat{I}_j(\lambda,\ha)=I_j(c^{\alpha}\lambda^{\alpha}x_0(\ha T_0),c^{\beta}\lambda^{\beta}y_0(\ha T_0))=O(\lambda^{-\frac{n+1}{n+2}-\varepsilon}),\hs \lambda\rightarrow +\infty,$$
   $$\hat{J}_j(\lambda,\ha)=J_j(c^{\alpha}\lambda^{\alpha}x_0(\ha T_0),c^{\beta}\lambda^{\beta}y_0(\ha T_0))=O(\lambda^{-\frac{1}{n+2}-\varepsilon}),\hs \lambda\rightarrow +\infty;$$

   {\rm (ii)} \hs When $r=1,s=0$, we have
  $$ \frac{\partial \hat{I}_j(\lambda,\ha)}{{\partial \lambda}}=\frac{\partial I_j(x,y)}{{\partial x}}\cdot\frac{\partial x}{{\partial \lambda}}+\frac{\partial I_j(x,y)}{{\partial y}}\cdot\frac{\partial y}{{\partial \lambda}},$$
   $$ \frac{\partial \hat{J}_j(\lambda,\ha)}{{\partial \lambda}}=\frac{\partial J_j(x,y)}{{\partial x}}\cdot\frac{\partial x}{{\partial \lambda}}+\frac{\partial J_j(x,y)}{{\partial y}}\cdot\frac{\partial y}{{\partial \lambda}},$$
   from (4.1), (4.2) and (2.6), we have\\
   $$ \frac{\partial \hat{I}_j(\lambda,\ha)}{{\partial \lambda}}=O(\lambda^{-\frac{n+1}{n+2}-\varepsilon-1}),\hs \lambda\rightarrow +\infty,$$
   $$ \frac{\partial \hat{J}_j(\lambda,\ha)}{{\partial \lambda}}=O(\lambda^{-\frac{1}{n+2}-\varepsilon-1}),\hs \lambda\rightarrow +\infty;$$

   {\rm (iii)} \hs When $r=0,s=1$, we have\\
  $$ \frac{\partial \hat{I}_j(\lambda,\ha)}{{\partial \ha}}=\frac{\partial I_j(x,y)}{{\partial x}}\cdot\frac{\partial x}{{\partial \ha}}+\frac{\partial I_j(x,y)}{{\partial y}}\cdot\frac{\partial y}{{\partial \ha}},$$
   $$ \frac{\partial \hat{J}_j(\lambda,\ha)}{{\partial \ha}}=\frac{\partial J_j(x,y)}{{\partial x}}\cdot\frac{\partial x}{{\partial \ha}}+\frac{\partial J_j(x,y)}{{\partial y}}\cdot\frac{\partial y}{{\partial \ha}},$$
   from (4.1), (4.2) and (2.6), we have\\
   $$ \frac{\partial \hat{I}_j(\lambda,\ha)}{{\partial \ha}}=O(\lambda^{-\frac{n+1}{n+2}-\varepsilon}),\hs \lambda\rightarrow +\infty,$$
   $$ \frac{\partial \hat{J}_j(\lambda,\ha)}{{\partial \ha}}=O(\lambda^{-\frac{1}{n+2}-\varepsilon}),\hs \lambda\rightarrow +\infty.$$\\

   Then, in the same way, for any non-negative integers $r$ and $s$($r+s\leq5$), we can prove\\
 $$\frac{\partial^{r+s} \hat{I}_j(\lambda,\ha)}{{\partial \lambda^{r}}{\partial \ha^{s}}}=O(\lambda^{-\frac{n+1}{n+2}-\varepsilon-r}),\hs \lambda\rightarrow +\infty,$$
   $$\frac{\partial^{r+s} \hat{J}_j(\lambda,\ha)}{{\partial \lambda^{r}}{\partial \ha^{s}}}=O(\lambda^{-\frac{1}{n+2}-\varepsilon-r}),\hs \lambda\rightarrow +\infty.$$ \qed \\

\noindent{\bf Lemma 4.2.} {\it  For $j=1,2,\cdots,k$, set $\lambda(t_j^-)=\lambda$, $\Delta\lambda(t_j)=\Delta\lambda$, $\ha(t_j^-)=\ha$, $\Delta\ha(t_j)=\Delta\ha$, then for any non-negative integers $r$ and $s$($r+s\leq5$),\\
$$\frac{\partial^{r+s} \Delta\lambda}{{\partial \lambda^{r}}{\partial \ha^{s}}}=O(\lambda^{-\varepsilon-r}),\hs \lambda\rightarrow +\infty,$$
   $$\frac{\partial^{r+s} \Delta\ha}{{\partial \lambda^{r}}{\partial \ha^{s}}}=O(\lambda^{-1-\varepsilon-r}),\hs \lambda\rightarrow +\infty,$$\\
if the condition (i) of Theorem 1.1 holds.}\\

\noindent{\bf Proof.}  Set $\hat{I}_j(\lambda,\ha)=\hat{I}_j$, $\hat{J}_j(\lambda,\ha)=\hat{J}_j$, $\lambda(t_j^-)=\lambda$, $\Delta\lambda(t_j)=\Delta\lambda$, $\ha(t_j^-)=\ha$, $\Delta\ha(t_j)=\Delta\ha$, then for large $\lambda$, by (2.6) and (2.12),  we have\\
\begin{eqnarray*} \bigtriangleup \lambda&=&\frac{1}{c}\{[x(t_j^-)+I_j]^{2n+2}+(n+1)[y(t_j^-)+J_j]^2\}^{\frac{n+2}{2n+2}}
  -\lambda\\&=& \frac{1}{c}\{[c^\frac{1}{n+2}\lambda^\frac{1}{n+2}x_0(\ha T_0)+\hat{I}_j]^{2n+2}+(n+1)[c^\frac{n+1}{n+2}\lambda^\frac{n+1}{n+2}y_0(\ha T_0)+\hat{J}_j]^2\}^{\frac{n+2}{2n+2}}
  -\lambda \\&=&\frac{1}{c}\{c^\frac{2n+2}{n+2}\lambda^\frac{2n+2}{n+2}x_0^{2n+2}(\ha T_0)+\sum_{i=1}^{2n+2}C_{2n+2}^{i}\hat{I}_{j}^{i}[c^\frac{1}{n+2}\lambda^\frac{1}{n+2}x_0(\ha T_0)]^{2n+2-i}\\
  &&+(n+1)[c^\frac{2n+2}{n+2}\lambda^\frac{2n+2}{n+2}y_0^{2}(\ha T_0)+2\hat{J}_jc^\frac{n+1}{n+2}\lambda^\frac{n+1}{n+2}y_0(\ha T_0)+\hat{J}_{j}^{2}]\}^{\frac{n+2}{2n+2}}
  -\lambda \\&=&\frac{1}{c}\{c^\frac{2n+2}{n+2}\lambda^\frac{2n+2}{n+2}+\sum_{i=1}^{2n+2}C_{2n+2}^{i}\hat{I}_{j}^{i}[c^\frac{1}{n+2}\lambda^\frac{1}{n+2} x_0(\ha T_0)]^{2n+2-i}\\&&+(n+1)[2\hat{J}_jc^\frac{n+1}{n+2}\lambda^\frac{n+1}{n+2}y_0(\ha T_0)+\hat{J}_{j}^{2}]\}^{\frac{n+2}{2n+2}}
  -\lambda \\&=&\lambda\{1+\sum_{i=1}^{2n+2}C_{2n+2}^{i}\hat{I}_{j}^{i}c^\frac{-i}{n+2}\lambda^\frac{-i}{n+2}x_0^{2n+2-i}(\ha T_0)+(2n+2)\hat{J}_j c^\frac{-n-1}{n+2}\lambda^\frac{-n-1}{n+2}y_0(\ha T_0)\\
  &&+(n+1)c^\frac{-2n-2}{n+2}\lambda^\frac{-2n-2}{n+2}\hat{J} _{j}^{2}\}^{\frac{n+2}{2n+2}}
  -\lambda\\&=&\lambda\{1+\sum_{i=1}^{2n+2}f_i(\ha)\hat{I}_{j}^{i}\lambda^\frac{-i}{n+2}+g_1(\ha)\hat{J}_j\lambda^\frac{-n-1}{n+2}+(n+1)c^\frac{-2n-2}{n+2}\hat{J} _{j}^{2}\lambda^\frac{-2n-2}{n+2}\}^{\frac{n+2}{2n+2}}
  -\lambda.\end{eqnarray*}\\
  Where $f_i(\ha)=C_{2n+2}^{i}c^\frac{-i}{n+2}x_0^{2n+2-i}(\ha T_0)$ and $g_1(\ha)=(2n+2)c^\frac{-n-1}{n+2}y_0(\ha T_0)$.\\
  By Taylor's formula, we have\\
  \begin{eqnarray*} \bigtriangleup \lambda&=&\lambda\{1+\sum_{i=1}^{+\infty}\frac{\frac{n+2}{2n+2}\cdot(\frac{n+2}{2n+2}-1)\cdot\cdot\cdot(\frac{n+2}{2n+2}+1-i)}{i!}[\sum_{i=1}^{2n+2}f_i(\ha)\hat{I}_{j}^{i}\lambda^\frac{-i}{n+2}\\&&+g_1(\ha)\hat{J}_j\lambda^\frac{-n-1}{n+2}+(n+1)c^\frac{-2n-2}{n+2}\hat{J} _{j}^{2}\lambda^\frac{-2n-2}{n+2}]^i\}
  -\lambda\\&=&\lambda\cdot\sum_{i=1}^{+\infty}\frac{\frac{n+2}{2n+2}\cdot(\frac{n+2}{2n+2}-1)\cdot\cdot\cdot(\frac{n+2}{2n+2}+1-i)}{i!}[\sum_{i=1}^{2n+2}f_i(\ha)\hat{I}_{j}^{i}\lambda^\frac{-i}{n+2}\\&&+g_1(\ha)\hat{J}_j\lambda^\frac{-n-1}{n+2}+(n+1)c^\frac{-2n-2}{n+2}\hat{J} _{j}^{2}\lambda^\frac{-2n-2}{n+2}]^i.\end{eqnarray*}\\
from Lemma 4.1, it is easily seen that for any non-negative integers $r$ and $s$($r+s\leq5$)\\
$$\frac{\partial^{r+s} \Delta\lambda}{{\partial \lambda^{r}}{\partial \ha^{s}}}=O(\lambda^{-\varepsilon-r}),\hs \lambda\rightarrow +\infty.\eqno (4.3)$$\\

  From (2.3) and (2.6), we have
  $$\left\{\begin{array}{ll}
  \hat{I}_j=c^{\alpha}\cdot(\lambda+\Delta\lambda)^{\alpha}\cdot x_0((\ha+\Delta\ha) T_0)-c^{\alpha}\cdot\lambda^{\alpha}\cdot x_0(\ha T_0),\\
  \hat{J}_j=c^{\beta}\cdot(\lambda+\Delta\lambda)^{\beta}\cdot y_0((\ha+\Delta\ha) T_0)-c^{\beta}\cdot\lambda^{\beta}\cdot y_0(\ha T_0).
  \end{array}\right.  \eqno (4.4)$$

     {\rm (i)} \hs When $|x_0(\ha T_0)|\leq\sqrt[2n+2]{\frac{1}{2}}$, $\sqrt{\frac{1}{2n+2}}\leq |y_0(\ha T_0)|\leq \sqrt{\frac{1}{n+1}}$, from (4.4), by Taylor's formula, we have\\
  $$\hat{I}_j=\alpha c^{\alpha}\cdot (\lambda+\xi_1\Delta\lambda)^{\alpha-1}\cdot x_0((\ha+\xi_1\Delta\ha)T_0)\cdot\Delta\lambda+T_0 c^{\alpha}\cdot (\lambda+\xi_1\Delta\lambda)^{\alpha}\cdot y_0((\ha+\xi_1\Delta\ha)T_0)\cdot\Delta\ha,\eqno (4.5)$$\\
   or
  $$T_0 c^{\alpha}\cdot (\lambda+\xi_1\Delta\lambda)^{\alpha}\cdot y_0((\ha+\xi_1\Delta\ha)T_0)\cdot\Delta\ha=\hat{I}_j-\alpha c^{\alpha}\cdot (\lambda+\xi_1\Delta\lambda)^{\alpha-1}\cdot x_0((\ha+\xi_1\Delta\ha)T_0)\cdot\Delta\lambda,$$\\
  where $0<\xi_1<1$,
  from Lemma 4.1 and (4.3), we have
  $$\hat{I}_j=O(\lambda^{\frac{-n-1}{n+2}-\varepsilon}), \hs \lambda\rightarrow +\infty,$$
  $$\alpha c^{\alpha}\cdot (\lambda+\xi_1\Delta\lambda)^{\alpha-1}\cdot x_0((\ha+\xi_1\Delta\ha)T_0)\cdot\Delta\lambda=O(\lambda^{\frac{-n-1}{n+2}-\varepsilon}), \hs \lambda\rightarrow +\infty.$$\\
  It is easily seen that $|y_0((\ha+\xi_1\Delta\ha)T_0)|>C$(C is a positive number), hence,\\
  $$\Delta\ha=O(\lambda^{-1-\varepsilon}), \hs \lambda\rightarrow +\infty.\eqno (4.6)$$\\
  From (4.5), we have
\begin{eqnarray*}\frac{\partial \hat{I}_j}{\partial \ha}&=&\alpha c^{\alpha}\cdot [\xi_1\cdot(\alpha-1)\cdot(\lambda+\xi_1\Delta\lambda)^{\alpha-2}\cdot\frac{\partial \Delta\lambda}{\partial \ha}]\cdot x_0((\ha+\xi_1\Delta\ha)T_0)\cdot\Delta\lambda\\&&+\alpha c^{\alpha}\cdot (\lambda+\xi_1\Delta\lambda)^{\alpha-1}\cdot[T_0\cdot y_0((\ha+\xi_1\Delta\ha)T_0)\cdot(1+\xi_1\frac{\partial \Delta\ha}{\partial \ha})]\cdot\Delta\lambda\\&&+\alpha c^{\alpha}\cdot (\lambda+\xi_1\Delta\lambda)^{\alpha-1}\cdot x_0((\ha+\xi_1\Delta\ha)T_0)\cdot\frac{\partial \Delta\lambda}{\partial \ha}\\&&+T_0 c^{\alpha}\cdot [\alpha \xi_1\cdot(\lambda+\xi_1\Delta\lambda)^{\alpha-1}\cdot\frac{\partial \Delta\lambda}{\partial \ha}]\cdot y_0((\ha+\xi_1\Delta\ha)T_0)\cdot\Delta\ha\\&&+T_0 c^{\alpha}\cdot (\lambda+\xi_1\Delta\lambda)^{\alpha}\cdot[-T_0\cdot x_0^{2n+1}((\ha+\xi_1\Delta\ha)T_0)\cdot(1+\xi_1\frac{\partial \Delta\ha}{\partial \ha})]\cdot\Delta\ha\\&&+T_0 c^{\alpha}\cdot (\lambda+\xi_1\Delta\lambda)^{\alpha}\cdot y_0((\ha+\xi_1\Delta\ha)T_0)\cdot\frac{\partial \Delta\ha}{\partial \ha},\end{eqnarray*}
i.e.,
\begin{eqnarray*}\frac{\partial \Delta\ha}{\partial \ha}&=&[T_0 c^{\alpha}\cdot (\lambda+\xi_1\Delta\lambda)^{\alpha}\cdot y_0((\ha+\xi_1\Delta\ha)T_0)\\&&+\alpha c^{\alpha}\cdot T_0\cdot\xi_1\cdot (\lambda+\xi_1\Delta\lambda)^{\alpha-1}\cdot y_0((\ha+\xi_1\Delta\ha)T_0)\cdot\Delta\lambda\\&&-T_0^2\cdot\xi_1\cdot c^{\alpha}\cdot (\lambda+\xi_1\Delta\lambda)^{\alpha}\cdot x_0^{2n+1}((\ha+\xi_1\Delta\ha)T_0)\cdot\Delta\ha]^{-1} \cdot\{\frac{\partial \hat{I}_j}{\partial \ha}\\&&-\alpha c^{\alpha}\cdot [\xi_1\cdot(\alpha-1)\cdot(\lambda+\xi_1\Delta\lambda)^{\alpha-2}\cdot\frac{\partial \Delta\lambda}{\partial \ha}]\cdot x_0((\ha+\xi_1\Delta\ha)T_0)\cdot\Delta\lambda\\&&-T_0\cdot\alpha c^{\alpha}\cdot (\lambda+\xi_1\Delta\lambda)^{\alpha-1}\cdot y_0((\ha+\xi_1\Delta\ha)T_0)\cdot\Delta\lambda\\&&-\alpha c^{\alpha}\cdot (\lambda+\xi_1\Delta\lambda)^{\alpha-1}\cdot x_0((\ha+\xi_1\Delta\ha)T_0)\cdot\frac{\partial \Delta\lambda}{\partial \ha}\\&&-T_0 c^{\alpha}\cdot [\alpha \xi_1\cdot(\lambda+\xi_1\Delta\lambda)^{\alpha-1}\cdot\frac{\partial \Delta\lambda}{\partial \ha}]\cdot y_0((\ha+\xi_1\Delta\ha)T_0)\cdot\Delta\ha\\&&+T_0^2 c^{\alpha}\cdot (\lambda+\xi_1\Delta\lambda)^{\alpha}\cdot x_0^{2n+1}((\ha+\xi_1\Delta\ha)T_0)\cdot\Delta\ha\},\end{eqnarray*}
from (4.3), (4.6) and Lemma 4.1, we have
$$\frac{\partial \Delta\ha}{\partial \ha}=O(\lambda^{-1-\varepsilon}), \hs \lambda\rightarrow +\infty.\eqno (4.7)$$

  In the same way, from (4.5), we have
\begin{eqnarray*}\frac{\partial \hat{I}_j}{\partial \lambda}&=&\alpha c^{\alpha}\cdot [(\alpha-1)\cdot(\lambda+\xi_1\Delta\lambda)^{\alpha-2}\cdot(1+\xi_1\frac{\partial \Delta\lambda}{\partial \lambda})]\cdot x_0((\ha+\xi_1\Delta\ha)T_0)\cdot\Delta\lambda\\&&+\alpha c^{\alpha}\cdot (\lambda+\xi_1\Delta\lambda)^{\alpha-1}\cdot[\xi_1\cdot T_0\cdot y_0((\ha+\xi_1\Delta\ha)T_0)\cdot\frac{\partial \Delta\ha}{\partial \lambda}]\cdot\Delta\lambda\\&&+\alpha c^{\alpha}\cdot (\lambda+\xi_1\Delta\lambda)^{\alpha-1}\cdot x_0((\ha+\xi_1\Delta\ha)T_0)\cdot\frac{\partial \Delta\lambda}{\partial \lambda}\\&&+T_0 c^{\alpha}\cdot [\alpha \cdot(\lambda+\xi_1\Delta\lambda)^{\alpha-1}\cdot(1+\xi_1\frac{\partial \Delta\lambda}{\partial \lambda})]\cdot y_0((\ha+\xi_1\Delta\ha)T_0)\cdot\Delta\ha\\&&+T_0 c^{\alpha}\cdot (\lambda+\xi_1\Delta\lambda)^{\alpha}\cdot[-T_0\cdot\xi_1\cdot x_0^{2n+1}((\ha+\xi_1\Delta\ha)T_0)\cdot\frac{\partial \Delta\ha}{\partial \lambda}]\cdot\Delta\ha\\&&+T_0 c^{\alpha}\cdot (\lambda+\xi_1\Delta\lambda)^{\alpha}\cdot y_0((\ha+\xi_1\Delta\ha)T_0)\cdot\frac{\partial \Delta\ha}{\partial \lambda},\end{eqnarray*}
i.e.,
\begin{eqnarray*}\frac{\partial \Delta\ha}{\partial \lambda}&=&[T_0 c^{\alpha}\cdot (\lambda+\xi_1\Delta\lambda)^{\alpha}\cdot y_0((\ha+\xi_1\Delta\ha)T_0)\\&&+\alpha c^{\alpha}\cdot T_0\cdot\xi_1\cdot (\lambda+\xi_1\Delta\lambda)^{\alpha-1}\cdot y_0((\ha+\xi_1\Delta\ha)T_0)\cdot\Delta\lambda\\&&-T_0^2\cdot\xi_1\cdot c^{\alpha}\cdot (\lambda+\xi_1\Delta\lambda)^{\alpha}\cdot x_0^{2n+1}((\ha+\xi_1\Delta\ha)T_0)\cdot\Delta\ha]^{-1} \cdot\{\frac{\partial \hat{I}_j}{\partial \lambda}\\&&-\alpha c^{\alpha}\cdot [(\alpha-1)\cdot(\lambda+\xi_1\Delta\lambda)^{\alpha-2}\cdot(1+\xi_1\frac{\partial \Delta\lambda}{\partial \lambda})]\cdot x_0((\ha+\xi_1\Delta\ha)T_0)\cdot\Delta\lambda\\&&-\alpha c^{\alpha}\cdot (\lambda+\xi_1\Delta\lambda)^{\alpha-1}\cdot x_0((\ha+\xi_1\Delta\ha)T_0)\cdot\frac{\partial \Delta\lambda}{\partial \lambda}\\&&-T_0 c^{\alpha}\cdot [\alpha \cdot(\lambda+\xi_1\Delta\lambda)^{\alpha-1}\cdot(1+\xi_1\frac{\partial \Delta\lambda}{\partial \lambda})]\cdot y_0((\ha+\xi_1\Delta\ha)T_0)\cdot\Delta\ha\},\end{eqnarray*}
from (4.3), (4.6) and Lemma 4.1, we have
$$\frac{\partial \Delta\ha}{\partial \lambda}=O(\lambda^{-2-\varepsilon}), \hs \lambda\rightarrow +\infty.\eqno (4.8)$$

      In the same way, for any non-negative integers $r$ and $s$($r+s\leq5$), we can prove
         $$\frac{\partial^{r+s} \Delta\ha}{{\partial \lambda^{r}}{\partial \ha^{s}}}=O(\lambda^{-1-\varepsilon-r}),\hs \lambda\rightarrow +\infty.\eqno (4.9)$$

      {\rm (ii)} \hs When $\sqrt[2n+2]{\frac{1}{2}}<|x_0(\ha T_0)|\leq1$, $|y_0(\ha T_0)|<\sqrt{\frac{1}{2n+2}}$, from (4.4), by Taylor's formula, we have
$$\hat{J}_j=\beta c^{\beta}(\lambda+\xi_2\Delta\lambda)^{\beta-1}y_0((\ha+\xi_2\Delta\ha)T_0)\cdot\Delta\lambda-T_0 c^{\beta}(\lambda+\xi_2\Delta\lambda)^{\beta}x_0^{2n+1}((\ha+\xi_2\Delta\ha)T_0)\cdot\Delta\ha,$$
 where $0<\xi_2<1$. In the same way, for any non-negative integers $r$ and $s$($r+s\leq5$), we can prove
         $$\frac{\partial^{r+s} \Delta\ha}{{\partial \lambda^{r}}{\partial \ha^{s}}}=O(\lambda^{-1-\varepsilon-r}),\hs \lambda\rightarrow +\infty.\eqno (4.10)$$

         Then, from (4.3), (4.9) and (4.10), Lemma 4.2 is proved.\qed \\

\noindent{\bf Lemma 4.3.} {\it  For $j=1,2,\cdots,k$, set $\mu(t_j^-)=\mu$, $\Delta\mu(t_j)=\Delta\mu$, $\phi(t_j^-)=\phi$, $\Delta\phi(t_j)=\Delta\phi$, then for any non-negative integers $r$ and $s$($r+s\leq5$),
$$\frac{\partial^{r+s} \Delta\mu}{{\partial \mu^{r}}{\partial \phi^{s}}}=O(\mu^{-\varepsilon-r}),\hs \mu\rightarrow +\infty,$$
   $$\frac{\partial^{r+s} \Delta\phi}{{\partial \mu^{r}}{\partial \phi^{s}}}=O(\mu^{-1-\varepsilon-r}),\hs \mu\rightarrow +\infty,$$
if the condition (i) of Theorem 1.1 holds.}\\

\noindent{\bf Proof.} From Lemma 3.2, there is a symplectic diffeomorphism $\psi_1$ depending periodically on $t$ of the form\\

$$\left\{\begin{array}{ll}
  \lambda=\tilde{\mu}+u(\tilde{\mu},\tilde{\phi},t),\\
  \theta=\tilde{\phi}+v(\tilde{\mu},\tilde{\phi},t),
\end{array}\right.  \eqno (4.11)$$\\
with $u\in F(1-(a-b))$ and $v\in F(-(a-b))$. Set $\tilde{\mu}(t_j^-)=\tilde{\mu}$, $\Delta\tilde{\mu}(t_j)=\Delta\tilde{\mu}$, $\tilde{\phi}(t_j^-)=\tilde{\phi}$, $\Delta\tilde{\phi}(t_j)=\Delta\tilde{\phi}$. Then we have
 $$\left\{\begin{array}{ll}
  \Delta\lambda=\Delta\tilde{\mu}+u(\tilde{\mu}+\Delta\tilde{\mu},\tilde{\phi}+\Delta\tilde{\phi})-u(\tilde{\mu},\tilde{\phi}),\\
  \Delta\theta=\Delta\tilde{\phi}+v(\tilde{\mu}+\Delta\tilde{\mu},\tilde{\phi}+\Delta\tilde{\phi})-v(\tilde{\mu},\tilde{\phi}),
\end{array}\right.  \eqno (4.12)$$\\
by Taylor's formula, we have\\
 $$\left\{\begin{array}{ll}
 \Delta\lambda=\Delta\tilde{\mu}+u_{\tilde{\mu}}(\tilde{\mu}+\xi_3\Delta\tilde{\mu},\tilde{\phi}+\xi_3\Delta\tilde{\phi})\cdot\Delta\tilde{\mu}+u_{\tilde{\phi}}(\tilde{\mu}+\xi_3\Delta\tilde{\mu},\tilde{\phi}+\xi_3\Delta\tilde{\phi})\cdot\Delta\tilde{\phi},\\
  \Delta\theta=\Delta\tilde{\phi}+v_{\tilde{\mu}}(\tilde{\mu}+\xi_4\Delta\tilde{\mu},\tilde{\phi}+\xi_4\Delta\tilde{\phi})\cdot\Delta\tilde{\mu}+v_{\tilde{\phi}}(\tilde{\mu}+\xi_4\Delta\tilde{\mu},\tilde{\phi}+\xi_4\Delta\tilde{\phi})\cdot\Delta\tilde{\mu},
\end{array}\right.  \eqno (4.13)$$\\
where $0<\xi_3<1$ and $0<\xi_4<1$. Then we have(set $u_{\tilde{\mu}}(\tilde{\mu}+\xi_3\Delta\tilde{\mu},\tilde{\phi}+\xi_3\Delta\tilde{\phi})=\frac{\partial u }{\partial \tilde{\mu}}$, $u_{\tilde{\phi}}(\tilde{\mu}+\xi_3\Delta\tilde{\mu},\tilde{\phi}+\xi_3\Delta\tilde{\phi})=\frac{\partial u}{\partial \tilde{\phi}}$, $v_{\tilde{\mu}}(\tilde{\mu}+\xi_4\Delta\tilde{\mu},\tilde{\phi}+\xi_4\Delta\tilde{\phi})=\frac{\partial v}{\partial \tilde{\mu}}$, $v_{\tilde{\phi}}(\tilde{\mu}+\xi_4\Delta\tilde{\mu},\tilde{\phi}+\xi_4\Delta\tilde{\phi})=\frac{\partial v}{\partial \tilde{\phi}}$)
$$\left\{\begin{array}{ll}
 \Delta\tilde{\mu}=[(1+\frac{\partial u }{\partial \tilde{\mu}})(1+\frac{\partial v}{\partial \tilde{\phi}})-\frac{\partial u}{\partial \tilde{\phi}}\cdot\frac{\partial v}{\partial \tilde{\mu}}]^{-1}\cdot[(1+\frac{\partial v}{\partial \tilde{\phi}})\cdot\Delta\lambda-\frac{\partial u}{\partial \tilde{\phi}}\cdot\Delta\ha],\\
  \Delta\tilde{\phi}=[(1+\frac{\partial u }{\partial \tilde{\mu}})(1+\frac{\partial v}{\partial \tilde{\phi}})-\frac{\partial u}{\partial \tilde{\phi}}\cdot\frac{\partial v}{\partial \tilde{\mu}}]^{-1}\cdot[(1+\frac{\partial u}{\partial \tilde{\mu}})\cdot\Delta\ha-\frac{\partial v}{\partial \tilde{\mu}}\cdot\Delta\lambda].
\end{array}\right.  \eqno (4.14)$$\\
From (4.11) and Lemma 4.2, we have\\
 $$\left\{\begin{array}{ll}
  \Delta\lambda=O(\lambda^{-\varepsilon})=O(\tilde{\mu}^{-\varepsilon}),\hs \tilde{\mu}\rightarrow+\infty,\\
  \Delta\theta=O(\lambda^{-1-\varepsilon})=O(\tilde{\mu}^{-1-\varepsilon}),\hs \tilde{\mu}\rightarrow+\infty.
\end{array}\right.  \eqno (4.15)$$\\
From (4.11), (4.14) and (4.15), we have\\
 $$\left\{\begin{array}{ll}
  \Delta\tilde{\mu}=O(\tilde{\mu}^{-\varepsilon}),\hs \tilde{\mu}\rightarrow+\infty,\\
  \Delta\tilde{\phi}=O(\tilde{\mu}^{-1-\varepsilon}),\hs \tilde{\mu}\rightarrow+\infty.
\end{array}\right.  \eqno (4.16)$$

Take the partial derivative of (4.12) with $\tilde{\mu}$, in the same way, we can prove
$$\left\{\begin{array}{ll}
  \frac{\partial \Delta\tilde{\mu} }{\partial \tilde{\mu}}=O(\tilde{\mu}^{-1-\varepsilon}),\hs \tilde{\mu}\rightarrow+\infty,\\
  \frac{\partial \Delta\tilde{\phi} }{\partial \tilde{\mu}}=O(\tilde{\mu}^{-2-\varepsilon}),\hs \tilde{\mu}\rightarrow+\infty.
\end{array}\right.  \eqno (4.17)$$

Take the partial derivative of (4.12) with $\tilde{\phi}$, in the same way, we can prove
$$\left\{\begin{array}{ll}
  \frac{\partial \Delta\tilde{\mu} }{\partial \tilde{\phi}}=O(\tilde{\mu}^{-\varepsilon}),\hs \tilde{\mu}\rightarrow+\infty,\\
  \frac{\partial \Delta\tilde{\phi} }{\partial \tilde{\phi}}=O(\tilde{\mu}^{-1-\varepsilon}),\hs \tilde{\mu}\rightarrow+\infty.
\end{array}\right.  \eqno (4.18)$$

 Then, in the same way, for any non-negative integers $r$ and $s$($r+s\leq5$), we can prove
$$\frac{\partial^{r+s} \Delta\tilde{\mu}}{{\partial \tilde{\mu}^{r}}{\partial \tilde{\phi}^{s}}}=O(\mu^{-\varepsilon-r}),\hs \mu\rightarrow +\infty,$$
   $$\frac{\partial^{r+s} \Delta\tilde{\phi}}{{\partial \tilde{\mu}^{r}}{\partial \tilde{\phi}^{s}}}=O(\mu^{-1-\varepsilon-r}),\hs \mu\rightarrow +\infty.$$

    Then, under the symplectic diffeomorphism $\psi_{m}\circ\psi_{m-1}\circ\cdot\cdot\cdot\circ\psi_{1}$ of Lemma 3.3, in the same way, for any non-negative integers $r$ and $s$($r+s\leq5$), we can prove
$$\frac{\partial^{r+s} \Delta\mu}{{\partial \mu^{r}}{\partial \phi^{s}}}=O(\mu^{-\varepsilon-r}),\hs \mu\rightarrow +\infty,$$
   $$\frac{\partial^{r+s} \Delta\phi}{{\partial \mu^{r}}{\partial \phi^{s}}}=O(\mu^{-1-\varepsilon-r}),\hs \mu\rightarrow +\infty.$$  \qed\\

\noindent {\bf 5. Proof of Theorems 1.1 - 1.2} \vskip 0.3cm

\noindent{\bf Lemma 5.1.} {\it The time 1 map $\Phi^1$ of the flow $\Phi^t$ of the vectorfield $X_{\hat{H}}$ given by  (3.3) is of the form
 $$\left\{\begin{array}{ll}
  \phi_1=\phi+r(\mu)+f(\mu,\phi),\\
  \mu_1=\mu+g(\mu,\phi),
\end{array}\right.  \eqno (5.1)$$
with $r(\mu)=a\mu^{a-1}+\int_{0}^{1}\frac{\partial \hat{h}_1(\mu,s)}{\partial \mu}ds$. Moreover, if $\mu$ is sufficiently large, then for for any non-negative integers $r$ and $s$($r+s\leq5$),
$$|D_{\mu}^{r}D_{\phi}^{s}f(\mu,\bullet)|,  |D_{\mu}^{r}D_{\phi}^{s}g(\mu,\bullet)|\leq \mu^{-\varepsilon_2},\hs \varepsilon_2>0,$$
  and if $\mu$ is sufficiently large, then $\dot{r}(\mu)>0$.}\\

\noindent{\bf Proof.} Let $(\mu(t),\phi(t))=(\mu(t,\mu,\phi),\phi(t,\mu,\phi))$ be the solution of (3.3),
 $(\mu(0),\phi(0))=(\mu,\phi)$. Set
$$r(\mu,t)=ta\mu^{a-1}+\int_0^t\frac{\partial \hat{h}_1(\mu,s)}{\partial \mu}ds.$$
  Then, we have
 $$\left\{\begin{array}{ll}
  \phi(t)=\phi+r(\mu,t)+A(\mu,\phi,t),\\
  \mu(t)=\mu+B(\mu,\phi,t),
\end{array}\right.  $$
where
\begin{eqnarray*} A(\mu,\phi,t)&= &a(a-1)\int_0^t\int_0^1(\mu+\tau B)^{a-2}\cdot Bd\tau ds+\int_0^t\int_0^1(\frac{\partial^2 \hat{h}_1}{\partial \mu^2})(\mu+\tau B,s)\cdot Bd\tau ds\\ &&+\int_0^t(\frac{\partial \hat{h}_2}{\partial \mu})(\mu+B,\phi+r+A,s)ds+C(t),\end{eqnarray*}
$B(\mu,\phi,t)=-\int_0^t(\frac{\partial \hat{h}_2}{\partial \phi})(\mu+B,\phi+r+A,s)ds+D(t),$\\
and
$$
C(t)=
\left\{
\begin{array}{ll}
  0,\hs 0\leq t<t_1,\\
  \sum_{i=1}^j\Delta\phi(t_i),\hs t_j\leq t<t_{j+1},\hs j=1,2,\cdots,k-1,\\
  \sum_{i=1}^k\Delta\phi(t_i),\hs t_k\leq t\leq 1,
\end{array}\right.  $$
$$
D(t)=
\left\{
\begin{array}{ll}
  0,\hs 0\leq t<t_1,\\
  \sum_{i=1}^j\Delta\mu(t_i),\hs t_j\leq t<t_{j+1},\hs j=1,2,\cdots,k-1,\\
  \sum_{i=1}^k\Delta\mu(t_i),\hs t_k\leq t\leq 1.
\end{array}\right.  $$
One verifies easily that for $\mu\geq\mu_0$ these equations have an unique solution in the space $|A|, |B|\leq1$. And we have $f(\mu,\phi)=A(\mu,\phi,1),  g(\mu,\phi)=B(\mu,\phi,1)$, Moreover, if $\mu$ is sufficiently large, by (3.3) and Lemma 3.1, we have

$$a(a-1)\int_0^1\int_0^1(\mu+\tau B)^{a-2}\cdot Bd\tau ds\in F(a-2),$$
$$\int_0^1\int_0^1(\frac{\partial^2 \hat{h}_1}{\partial \mu^2})(\mu+\tau B,s)\cdot Bd\tau ds\in F(b-2),$$
$$\int_0^1(\frac{\partial \hat{h}_2}{\partial \mu})(\mu+B,\phi+r+A,s)ds\in F(-\varepsilon_1-1),$$
$$-\int_0^1(\frac{\partial \hat{h}_2}{\partial \phi})(\mu+B,\phi+r+A,s)ds\in F(-\varepsilon_1),$$
and by Lemma 4.3, for any non-negative integers $r$ and $s$($r+s\leq5$), we have
$$\frac{\partial^{r+s} D(1)}{{\partial \mu^{r}}{\partial \phi^{s}}}=O(\mu^{-\varepsilon-r}),$$
   $$\frac{\partial^{r+s} C(1)}{{\partial \mu^{r}}{\partial \phi^{s}}}=O(\mu^{-1-\varepsilon-r}).$$
Then, if $\mu$ is sufficiently large, for any non-negative integers $r$ and $s$($r+s\leq5$), we have
$$|D_{\mu}^{r}D_{\phi}^{s}f(\mu,\bullet)|,  |D_{\mu}^{r}D_{\phi}^{s}g(\mu,\bullet)|\leq \mu^{-\varepsilon_2},\hs 0<\varepsilon_2<\min(2-a,\varepsilon,\varepsilon_1).$$
  From
  $$ r(\mu)=a\mu^{a-1}+\int_{0}^{1}\frac{\partial \hat{h}_1(\mu,s)}{\partial \mu}ds, $$
  we have
  $$ \dot{r}(\mu)=a(a-1)\mu^{a-2}+\int_{0}^{1}\frac{\partial \hat{h}_1^2(\mu,s)}{\partial \mu^2}ds,$$
  by (3.3), we have $\int_{0}^{1}\frac{\partial \hat{h}_1^2(\mu,s)}{\partial \mu^2}ds\in F(b-2)$, and $\mu^{a-2}\in F(a-2)$, hence, if $\mu$ is sufficiently large, $\dot{r}(\mu)>0$.\\

\noindent{\bf Lemma 5.2.} {\it If the condition (ii) of Theorem 1.1 holds,
  the time 1 map $\Phi^1$ of (3.3) is area-preserving. Moreover, the map $\Phi^1$ has the intersection property on $A_{\mu_0}$, i.e. if C is an embedded circle in $A_{\mu_0}$ homotopic to a circle $\mu=const$ in $A_{\mu_0}$, then $\Phi(C)\bigcap C\neq\emptyset$.}\\

\noindent{\bf Proof.}   The time 1 map $\Phi^1$ of (3.3) is
$$\Phi^1=P_{k}\circ\Phi^{*}_{k}\circ\cdot\cdot\cdot\circ P_1\circ \Phi^{*}_{1}\circ P_0,$$
where
$$\Phi^{*}_{j}:(\mu(t_j^-),\phi(t_j^-))\mapsto (\mu(t_j^-)+\bigtriangleup \mu(t_j),\phi(t_j^-)+\bigtriangleup \phi(t_j)),\hs j=1,2,\cdots,k;$$
$$P_{k}:(\mu(t_{k}),\phi(t_{k}))\mapsto (\mu(1),\phi(1));$$
$$P_{j}:(\mu(t_j),\phi(t_j))\mapsto (\mu(t_{j+1}^-),\phi(t_{j+1}^-)),j=1,2,\cdots,k-1;$$
$$P_{0}:(\mu(0),\phi(0))\mapsto (\mu(t_1^-),\phi(t_1^-)).$$\\
 The condition (ii) of Theorem 1.1 implies the Jacobian $|\Delta|=1$ of jump maps $\Phi_{j}:(x,y)\mapsto (x,y)+(I_j(x,y),J_j(x,y)),\hs j=1,2,\cdots,k$. Thus $\Phi_{j}$ are area-preserving. Under the symplectic diffeomorphism $\psi_{m}\circ\psi_{m-1}\circ\cdot\cdot\cdot\circ\psi_{1}\circ\va_0$, $\Phi^{*}_{j}$ are also area-preserving; Since $P_{j}(j=0,1,2,\cdots,k)$ are the Poincar$\acute{e}$ maps of Hamiltonian equations (3.3) without  impulsive terms, then $P_{j}( j=0,1,2,\cdots,k)$ are area-preserving, hence, the map $\Phi^1$ is area-preserving, this implies that the map $\Phi^1$ has the intersection property. \\

\noindent{\bf Lemma 5.3} (Moser's theorem [14, 15]).{\it Consider an annulus
  $$D:a\leq \mu\leq b, \hs \phi\in S^1,$$
   where $b-a>0$. Let $P: D \mapsto D$ be a map defined by
  $$\left\{\begin{array}{ll}
  \mu_1=\mu+f(\mu,\phi),\\
  \phi_1=\phi+\alpha(\mu)+g(\mu,\phi),
  \end{array}\right.  $$\\
  where $f(\mu,\phi)$,$g(\mu,\phi)$ are 1-periodic functions with respect to $\phi$.\\
Assume:\\
(\romannumeral1). $\dot{\alpha}(\mu)>0$, $\forall a\leq \mu \leq b;$(twist condition)\\
(\romannumeral2). for $\forall (\mu,\phi)\in D$,
  $$\sum_{m+n\leq 5}| \frac{\partial^{m+n} f(\mu,\phi)}{{\partial \mu^m}{\partial \phi^n}}|, \sum_{m+n\leq 5}| \frac{\partial^{m+n} g(\mu,\phi)}{{\partial \mu^m}{\partial \phi^n}}| \leq \varepsilon_0,$$
   where $\varepsilon_0$ is a very small constant;(small property condition)\\
(\romannumeral3). The map $P$ has the intersection property, i.e. each closed curve $C$ in $D$ homotopic to a circle $\mu=1$, then
  $$P(C)\bigcap C \neq \emptyset.$$
Then, there exists $\bar{\mu}\in [a,b]$ and a closed curve $C$ such that\\
(\romannumeral1). $P(C)=C$, i.e. $C$ is an invariant curve of $P$,\\
(\romannumeral2). $P|c: C \mapsto C,$  is a rotation with the frequency $\omega=\alpha(\bar{\mu})$.}\\

We are now in a position to prove our main results Theorem 1.1 and 1.2 mentioned in the section 1.\\

Let
$$\Omega_m=\left\{(\phi,\mu)|\phi\in S^1, m\leq\mu\leq m+1\right\}$$
and let $\Phi^1_m$ be the restriction of $\Phi^1$ in $\Omega_m$.

   For any $m>m_0$($m_0$ is  large enough), by Lemma 5.1 and Lemma 5.2, $\Phi^1_m$ satisfies all conditions of Lemma 5.3. Then, there exists an invariant curve $C_m$ in $\Omega_m$ and the map $\Phi^1$ restricted on $C_m$ has a rotation with the frequency $\omega=\alpha(\bar{\mu})$($\bar{\mu}\in [m,m+1]$, $\alpha(\bar{\mu})>\alpha(m_0)$). Hence, there are many invariant curves with $m>m_0$.

   For any $\omega$, there is an embedding $\varphi$: $S^1\rightarrow A_{\mu_0}$ of a circle, which is differentiably close to the injection map $j$ of the circle $\left\{\omega\right\}\times S^1\rightarrow A_{\mu_0}$, and which is invariant under the map $F=\Phi^1$. Moreover, on this invariant curve the map $F$ is conjugated to a rotation with rotation number $\omega$:
   $$F\circ\varphi(s)=\varphi(s+\omega)\hs with \hs s(mod \hs 1).\eqno (5.2)$$
The solutions of the Hamitonian equation starting at time $t=0$ on this invariant curve determine a 1-periodic cylinder in the space $(\mu,\phi,t)\in A_{\mu_0}\times R$. Since the Hamiltonian vectorfield $X_{\hat{H}}$ is timeperiodic, the phase space is $A_{\mu_0}\times S^1$. Let $\Phi^t_1$ with $\Phi^0_1=id$ be the flow of the time-independent vectorfield $(X_{\hat{H}},1)$ on $A_{\mu_0}\times S^1$ and define the embedded torus $\Psi$: $T^2\rightarrow A_{\mu_0}\times S^1$ by setting
$$\Psi(s,\tau)=\Phi^\tau_1(\varphi(s-\tau\omega),0)=(\Phi^\tau\circ\varphi(s-\tau\omega),\tau).\eqno (5.3)$$
In view of (5.2) we have with $F=\Phi^1$ indeed $\Psi(s+1,\tau)=\Psi(1,\tau+1)=\Psi(s,\tau)$. Moreover $\Phi^t_1\circ\Psi(s,\tau)=\Psi(s+\omega t,\tau+t)$, so that the torus $\Psi(T^2)$ is quasiperiodic having the frequencies $(\omega,1)$. This proves the statement of Theorem 1.2. In order to prove the statement of Theorem 1.1 just observed that, in the original coordinates, every point $(x,y)\in R^2$ is in the interior of some invariant curve of the time 1 map of the flow which goes around the origin. Its solution is therefore confined in the interior of the time periodic cylinder above the invariant curve and hence is bounded. This ends the proof of Theorem 1.1.\\

\noindent{\bf 6. Further discussions about the jump maps $I_j, J_j$} \\

From Theorem 1.1 and 1.2, the question arise naturally about what kind of jump maps $I_j(x,y)$ and $J_j(x,y)$ can satisfy conditions (i) and (ii) of Theorem 1.1 and 1.2. In this section, we will give such jump maps as examples. We should mention that,
  it is hard to take examples for giving the jump maps $I_j(x,y)$ and $J_j(x,y)$  directly in $(x,y)$-plane, but under $\va_0$, it is not difficult. In fact, from Lemma 4.2 and the condition (i) of Theorems, for any non-negative integers $r$ and $s$($r+s\leq5$), we have
  $$\left\{\begin{array}{ll}
\frac{\partial^{r+s} \Delta\lambda}{{\partial \lambda^{r}}{\partial \ha^{s}}}=O(\lambda^{-\varepsilon-r}),\hs \lambda\rightarrow +\infty,\\
   \frac{\partial^{r+s} \Delta\ha}{{\partial \lambda^{r}}{\partial \ha^{s}}}=O(\lambda^{-1-\varepsilon-r}),\hs \lambda\rightarrow +\infty,
 \end{array}\right.  \eqno (6.1)$$
 the condition (ii) of Theorems 1.1-1.2 is the same as jump maps $$  \Phi_j:(x,y)\mapsto (x,y)+(I_j(x,y),J_j(x,y)),$$ are area-preserving, and $\va_0$ is a symplectic diffeomorphism, so jump maps $$ (\lambda,\ha)\mapsto (\lambda,\ha)+(\Delta\lambda,\Delta\ha)$$ are also area-preserving. It is the same as
 $$\frac{\partial \bigtriangleup\lambda}{\partial \lambda}+\frac{\partial \bigtriangleup\ha}{\partial \ha}+\frac{\partial \bigtriangleup\lambda}{\partial \lambda}\cdot \frac{\partial \bigtriangleup\ha}{\partial \ha}-\frac{\partial \bigtriangleup\lambda}{\partial \ha}\cdot \frac{\partial \bigtriangleup\ha}{\partial \lambda}=0.\eqno (6.2)$$\\
 We can list some examples easily. For example:\\
     $$\left\{\begin{array}{ll}
  \bigtriangleup \ha(t_j)=\frac{1}{\lambda^2(t_j^-)},\\
  \bigtriangleup \lambda(t_j)=0.
  \end{array}\right.  \eqno (6.3)$$\\

  It is easy to know (6.1) and (6.2) are satisfied with impulse (6.3). More examples can also given, for example:\\

\noindent (i) $\bigtriangleup \lambda(t_j)=0,\hs  \bigtriangleup \ha(t_j)=\frac{1}{\lambda^3(t_j^-)};$\\
\noindent (ii) $\bigtriangleup \lambda(t_j)=0,\hs  \bigtriangleup \ha(t_j)=\frac{1}{\lambda^4(t_j^-)}$, etc.\\

  Now we prove the existence of the impulse in $(x,y)$-plane for the case (6.3). \\

  We set up $I_j(x(t_j^-),y(t_j^-))=I_j$, $J_j(x(t_j^-),y(t_j^-))=J_j$, $\lambda(t_j)=\lambda_j$, $\ha(t_j)=\ha_j$, $x(t_j)=x_j$, $y(t_j)=y_j$, $\lambda(t_{j^-})=\lambda_{j^-}$, $\ha(t_{j^-})=\ha_{j^-}$, $x(t_{j^-})=x_{j^-}$, $y(t_{j^-})=y_{j^-}$. Under $\va_0$, we know $(\lambda_{j^-},\ha_{j^-})\mapsto (x_{j^-},y_{j^-})$ and $(\lambda_{j},\ha_{j})\mapsto (x_{j},y_{j})$. We set up map Q: $(\lambda_{j^-},\ha_{j^-})\mapsto (\lambda_{j},\ha_{j})$. It is easy to know that any map $\Phi:(x_{j^-},y_{j^-})\mapsto (x_{j},y_{j})$ has:
  $$\Phi=\va_0\circ Q\circ\va_0^{-1}.$$
  For $\forall (x_{j^-},y_{j^-})$, we have
  \begin{eqnarray*} \Phi(x_{j^-},y_{j^-})&= & \va_0\circ Q\circ\va_0^{-1}(x_{j^-},y_{j^-})\\&= &\va_0\circ Q(\lambda_{j^-},\ha_{j^-})\\
&=&\va_0(\lambda_{j},\ha_{j}).\end{eqnarray*}
  So the existence is proved. \\

  Now we calculate the jump map $I_j$ and $J_j$ for the case (6.3).

  From (2.6), the Jacobian $|\Delta|\neq 0$, by means of the implicit function theorem, and by (2.9), we can get\\
$$\ha=T_0^{-1}\cdot x_0^{-1}(c^\frac{-1}{n+2}\lambda^\frac{-1}{n+2}x)=T_0^{-1}\cdot x_0^{-1}(\frac{x}{[x^{2n+2}+(n+1)y^2]^{\frac{1}{2n+2}}}).\eqno (6.4)$$
Then from (2.9) and (6.4), we have
  $$\left\{\begin{array}{ll}
  \lambda_{j^-}=\frac{[x_{j^-}^{2n+2}+(n+1)y_{j^-}^2]^{\frac{n+2}{2n+2}}}{c},\\
  \ha_{j^-}=T_0^{-1}\cdot x_0^{-1}(\frac{x_{j^-}}{[x_{j^-}^{2n+2}+(n+1)y_{j^-}^2]^{\frac{1}{2n+2}}}).
  \end{array}\right.  \eqno (6.5)$$

  From (6.3) and (6.5), we have
  $$\left\{\begin{array}{ll}
  \lambda_{j}=\frac{[x_{j^-}^{2n+2}+(n+1)y_{j^-}^2]^{\frac{n+2}{2n+2}}}{c},\\
  \ha_{j}=T_0^{-1}\cdot x_0^{-1}(\frac{x_{j^-}}{[x_{j^-}^{2n+2}+(n+1)y_{j^-}^2]^{\frac{1}{2n+2}}})+\frac{c^2}{[x_{j^-}^{2n+2}+(n+1)y_{j^-}^2]^{\frac{n+2}{n+1}}}.
  \end{array}\right.  \eqno (6.6)$$

  From (2.6) and (6.6), we have
  $$\left\{\begin{array}{ll}
  x_{j}=[x_{j^-}^{2n+2}+(n+1)y_{j^-}^2]^{\frac{1}{2n+2}}\cdot x_0(x_0^{-1}(\frac{x_{j^-}}{[x_{j^-}^{2n+2}+(n+1)y_{j^-}^2]^{\frac{1}{2n+2}}})+\frac{T_0c^2}{[x_{j^-}^{2n+2}+(n+1)y_{j^-}^2]^{\frac{n+2}{n+1}}}),\\
  y_{j}=[x_{j^-}^{2n+2}+(n+1)y_{j^-}^2]^{\frac{1}{2}}\cdot y_0(x_0^{-1}(\frac{x_{j^-}}{[x_{j^-}^{2n+2}+(n+1)y_{j^-}^2]^{\frac{1}{2n+2}}})+\frac{T_0c^2}{[x_{j^-}^{2n+2}+(n+1)y_{j^-}^2]^{\frac{n+2}{n+1}}}).
\end{array}\right.  \eqno (6.7)$$\\

From (6.7), we have the jump map $I_j$ and $J_j$:
 \begin{eqnarray*}I_j&=&[x_{j^-}^{2n+2}+(n+1)y_{j^-}^2]^{\frac{1}{2n+2}}\cdot x_0(x_0^{-1}(\frac{x_{j^-}}{[x_{j^-}^{2n+2}+(n+1)y_{j^-}^2]^{\frac{1}{2n+2}}})\\
 &&+\frac{T_0c^2}{[x_{j^-}^{2n+2}+(n+1)y_{j^-}^2]^{\frac{n+2}{n+1}}})-x_{j^-},\end{eqnarray*}
 \begin{eqnarray*}J_j&=&[x_{j^-}^{2n+2}+(n+1)y_{j^-}^2]^{\frac{1}{2}}\cdot y_0(x_0^{-1}(\frac{x_{j^-}}{[x_{j^-}^{2n+2}+(n+1)y_{j^-}^2]^{\frac{1}{2n+2}}})\\
 &&+\frac{T_0c^2}{[x_{j^-}^{2n+2}+(n+1)y_{j^-}^2]^{\frac{n+2}{n+1}}})-y_{j^-}.\end{eqnarray*}\\

\vskip 0.5cm
\noindent{\bf Acknowledgments}\\

 This work is supported by the National Natural Science Foundation of China (No. 11571088) and the Zhejiang Provincial Natural Science Foundation of China (LY14A010024).\\

\noindent{\bf References} \\
  \newcounter{cankao}
\begin{list}
{[\arabic{cankao}]}{\usecounter{cankao}\itemsep=0cm} \small
\item M. Akhmet, Principles of Discontinuous Dynamical Systems, Springer-New York, 2010.
\item V. Alekseev, Quasirandom dynamical systems, 1, 2, 3, Math. USSR-Sb., 5 (1968) 73-128; 6 (1968) 505-560; 7 (1969) 1-43.
\item V. Arnold, Mathematical Methods of Classical Mechanics, Higher Education Press, Beijing, 2006.
\item
R. Dieckerhoff and E. Zehnder, Boundedness of solutions via the Twist
Theorem, Ann. Sc. Norm.Super. Pisa 14(1)(1987), 79-85.
 \item Y. Dong, Sublinear impulse effects and solvability of boundary value problems for differential equations with impulses, J. Math. Anal. Appl. 264 (2001), 32-48.
\item F. Jiang, J. Shen and Y. Zeng, Applications of the Poincar谷每Birkhoff theorem to impulsive Duffing equations at reso-nance, Nonlinear Anal. Real World Appl. 13 (2012), 1292每1305.
\item S. Laederich and M. Levi, Invariant curves and time-dependent potential, Ergodic. Theory Dynam. Systems 11 (1991), 365-378.
     \item V. Lakshmikantham, D. Bainov and P. Simeonov, Theory of Impulsive Differential Equations, World Scientific, Singapore, 1989.
  \item M. Levi, Quasiperiodic motions in superquadratic time periodic potentials. Comm. Math. Phys. 143(1) (1991), 43-83.

    \item J. Littlewood, Some problems in real and complex analysis, Heath, Lexington, Mass. 1968.
    \item B. Liu, Boundedness for solutions of nonlinear Hill's equations with periodic forcing terms via Moser's twist theorem, Journal of Differential Equations 79(1989), 304-315.
\item B. Liu, Boundedness of Solutions of Nonlinear Differential Equations, Journal of Differential Equations 144 (1998), 66-98.

  \item G. Morris, A case of boundedness of Littlewood's problem on oscillatory differential equations, Bull.
Austral. Math. Soc. 14 (1976), 71-93.

\item J. Moser, On invariant curves of aera-preserving mapping of annulus, Nachr. Akad. Wiss. Gottingen
Math. Phys. 2 (1962), 1-20.
\item J. Moser, Stable and Random Motion in Dynamic Systems, Ann. of Math. Studies, Princeton Uni. Press,
Princeton, NJ, 1973.
\item J. Norris, Boundedness in periodically forced second order conservative systems, J. London Math. Soc. 45 (2) (1992), 97-112.
 \item J. Nieto, Basic theory for nonresonace impulsive periodic problems of first order, J. Math. Anal. Appl. 205 (1997) 423每433.
  \item J. Nieto and D. O'Regan, Variational approach to impulsive differential equations, Nonlinear Anal. Real World Appl. 10 (2009), 680每690.
  \item R. Ortega, Asymmetric oscillators and twist mappings, J. Lond. Math. Soc. 53 (1996) 325-342.
 \item D. Qian, L. Chen and X. Sun, Periodic solutions of superlinear impulsive differential equations: A geometric approach, J. Differential Equations 258 (2015), 3088-3106.
 \item H. Russman, Uber invariant Kurven differenzierbarer Abbildungen eines Kreisringes, Nachr. Akad. Wiss. Gottingen Math. Phys. 2(1970), 67-105.
\item K. Sitnikov, Existence of oscillating motions for three-body problem, Dokl. Akad. Nauk., SSSR., 133 (2) (1960), 303-306.
\item J. Sun, H. Chen and J. Nieto, Infinitely many solutions for second-order Hamiltonian system with impulsive effects, Math. Comput. Modelling 54 (2011), 544每555

 \item X. Yuan, Lagrange stability for Duffing-type equations, J. Differential Equations 160 (1) (2000) 94-117.
\item X. Yuan, Invariant Tori of Duffing-type equations, J. Diff
erential Equations 142 (2) (1998), 231-262.
\item R. Yuan and X. Yuan, Boundedness of solutions for a class of nonlinear dfferential equations of second order via Moser＊s twist theorem, Nonlinear Analysis 46 (2001), 1073-1087.
\item F. Zhang, Quasiperiodic solutions of higher dimensional Duffing＊s equations via the KAM theorem, Sci. China Ser. A 44 (5) (2001), 631-644.

\item X. Yuan, Boundedness of solutions for Duffing equation with low regularity
in time, http://arxiv.org/abs/1705.02881v1, May, 2017

\end{list}
\end{document}